\documentclass[12pt]{article}
\usepackage{amsfonts}
\usepackage{amssymb,a4}
\usepackage{bbm}
\usepackage{amsmath,amsfonts,amssymb,amsthm}
\newcommand{\al}{\alpha}
\newcommand{\ot}{\otimes}

\def\wt{{\rm wt}}

\def\de{\delta}
\def\be{\beta}

\newcommand{\la}{\lambda}

\def\C{{\mathbb C}}

\def\Z{{\mathbb Z}}

\def\N{{\mathbb N}}
\def\1{{\bf 1}}

\def \wt{{\rm wt}}

\def \End{{\rm End}}

\def \<{\langle}
\def \>{\rangle}
\def \w{\omega}
\def \pf{\noindent {\bf Proof: \,}}

\def\ot{\otimes}
\def\om{\omega}

\def\al{\alpha}

\def\be{\beta}
\def \bee{\begin{equation}\label}

\def\de{\delta}
\def\pa{\partial}

\def\la{\lambda}

\def\k{\kappa}

\def\d{\mathbbm{d}}
\def\cc{\mathbbm{c}}

\def\qed{\hfill\mbox{$\Box$}}

\def\bZ{{\mathbb Z}}

\def\bC{{\mathbb C}}
\def\cc{\mathbbm{c}}

\def\C{{\mathscr C}}

\def\b{\mathfrak{b}}

\def\Ind{\mbox{\rm Ind}\,}

\def\End{\text{\rm End}}

\def\wt{\text{\rm wt}}

\def\C{{\mathbb C}}

\def\Z{{\mathbb Z}}

\def\N{{\mathbb N}}
\def\1{{\bf 1}}
\def\b1{\mathbbold{1}}
\def\k{{\bf k}}
\def \End{{\rm End}}

\def \bla{\tilde{\lambda}}
\def \bm{\tilde{\mu}}
\def \bn{\tilde{\nu}}
\def\H{\widehat{H}_{4}}
\def\NO{\mbox{\,$\circ\atop\circ$}\,}

\newtheorem{theorem}{Theorem}[section]
\newtheorem{prop}[theorem]{Proposition}
\newtheorem{lem}[theorem]{Lemma}
\newtheorem{coro}[theorem]{Corollary}
\newtheorem{rem}[theorem]{Remark}
\theoremstyle{definition}
\newtheorem{definition}[theorem]{Definition}

\numberwithin{equation}{section}

\begin{document}
\begin{center}
{\Large {\bf Extension of vertex operator algebra
$V_{\widehat{H}_{4}}(\ell,0)$}} \\
\vspace{0.5cm}  Cuipo Jiang\footnote{Supported  by China NSF grants
10931006,10871125, the Innovation Program of Shanghai Municipal
Education Commission (11ZZ18) and a grant of Science and Technology
Commission
of Shanghai Municipality (No. 09XD1402500)(E-mail: cpjiang@sjtu.edu.cn).} and Song Wang\footnote{E-mail: wangsong025@sina.com.}\\
 Department of Mathematics, Shanghai Jiaotong University, Shanghai 200240 China
\end{center}
\hspace{1cm}

\begin{abstract}
We classify  the irreducible restricted modules for the affine
Nappi-Witten Lie algebra $\H$ with some natural conditions. It turns
out the representation theory of $\H$ is quite different from the
theory of representations of Heisenberg algebras. We also study the
extension
 of the vertex operator
algebra $V_{\widehat{H}_{4}}(\ell,0)$ by the even lattice $L$. We
give  the structure of  the extension
$V_{\widehat{H}_{4}}(\ell,0)\otimes \C[L]$ and its irreducible
modules via irreducible representations of
$V_{\widehat{H}_{4}}(\ell,0)$ viewed as a vertex algebra.

\vskip 0.3cm
 2000MSC:17B69

\end{abstract}

\section{Introduction}

The Nappi-Witten algebra $H_{4}$ is a four-dimensional Lie algebra
with the basis $\{ a, b,c,d\}$ and the following Lie bracket:
$$
[a,b]=c,\ \ [d,a]=a,\ \ [d,b]=-b,\ \ [c,d]=[c,a]=[c,b]=0.
$$
There exists a non-degenerate symmetric bilinear form $(,)$ on
$H_{4}$ defined by
$$
(a,b)=1,\quad (c,d)=1,\quad\text{otherwise}, \ (\ ,\ )=0.
$$
 The
associated  affine Nappi-Witten Lie algebra  is
$$
\widehat{H}_{4}=H_{4}\ot\bC[t,t^{-1}]\oplus\bC{\bf k}.
$$

 The study of the
affine Nappi-Witten algebra $\H$ was initiated by the theory of
Wess-Zumino-Novikov-Witten (WZNW) models  based on non-abelian
non-semisimple groups \cite{Witten, NW}. WZNW models corresponding
to abelian or compact semisimple groups have been studied
extensively.  The importance of non-abelian non-semisimple types was
first noticed in \cite{Di, NW}. It was shown in \cite{NW} that the
NW model corresponding to a central extension of the two-dimensional
Euclidean groups describes the spectrum and scattering of  strings
moving in a gravitational wave background. In \cite {KK},
irreducible representations with unitary base for the affine NW
algebra were constructed through four free fields. In \cite{BJP},
Verma type modules for the affine Nappi-Witten algebra $\H$ were
investigated.
 The associated vertex operator algebra $V_{\H}(\ell,0)$ for a given complex number $\ell$
 and
its representations were also studied in \cite{BJP}. More results on
the NW model can be found in \cite{CFS,DQ1,DQ2,DK}.

It is known that every restricted $\H$-module of level $\ell$
carries a unique module structure for the vertex algebra
$V_{\H}(\ell,0)$, extending the action of $\H$ in a canonical way
and that every module of the vertex algebra $V_{\H}(\ell,0)$ is
naturally a restricted $\H$-module of level $\ell$. This leads us to
study restricted modules for the affine Lie algebra $\H$. On the
other hand, let  $L=\Z c\otimes \Z d$, then $L$ is a non-degenerate
even lattice. Motivated by \cite{Li2}, \cite{Li3}, \cite{Li4},
\cite{LW} and \cite{DLM1}, we consider the extension
$V_{\H}(\ell,0)[L]$ of $V_{\H}(\ell,0)$ by $L$. Some ideas we use
come from these papers. For the simple vertex operator algebra
$V_{\hat{\frak g}}(\ell,0)$ and an irreducible highest weight
$V_{\hat{\frak g}}(\ell,0)$-module $W$ with ${\frak g}$ a
finite-dimensional simple Lie algebra and $\ell$ a positive integer,
or with ${\frak g}$ abelian and $\ell=1$,  it is known that
$W^{(\al)}$ is still an irreducible highest weight module for the
affine Lie algebra $\hat{\frak g}$ (cf. \cite{DLM1}, \cite{Li2},
\cite{Li3}, \cite{Li4}, \cite{LW}), where $\al\in {\frak h}$ such
that $\al(0)$ acts semisimply on $V$ with integral eigenvalues and
${\frak h}$ is a Cartan subalgebra of ${\frak g}$. But for
$V_{\H}(\ell,0)$, it turns out it is quite different. For general
$\al\in L$ and an irreducible highest weight $\H$-module $W$,
$W^{(\al)}$ is not a highest weight module of $\H$ anymore. We
characterize $W^{(\al)}$ for every irreducible highest weight module
of $\H$.

The paper is organized as follows: In Section 2, we recall the
various modules, admissible modules and ordinary modules for a
vertex operator algebra. We also give some  known related results on
the affine Nappi-Witten algebra and the associated vertex operator
algebra. In section 3, we characterize three types of restricted
$\H$-modules and their irreducible quotients. Section 4 is devoted
to the extension $V_{\H}(\ell,0)[L]$ and its irreducible modules.
 \vspace{0.5cm}
\section{Preliminaries}

In this section, we briefly review the definitions of vertex
operator algebra and various modules (cf. \cite{Bo}, \cite{D},
\cite{DLM2}, \cite{FLM}, \cite{Z}, \cite{LL}, \cite{Xu}). We also
give some known results related to the Nappi-Witten algebra $H_{4}$.

A  vertex operator algebra $V$ is a $\Z$-graded vector space
$V=\bigoplus_{n\in \Z} V_n$ equipped with a linear map $Y:V\to(\End
V)[[z,z^{-1}]],\,a\mapsto Y(a,z)=\sum\limits_{n\in \Z}a_nz^{-n-1}$
for $a\in V$ such that $\dim V_n$ is finite for all integer $n$ and
that $V_n=0$ for sufficiently small  $n$ (see \cite{FLM}). There are
two distinguished vectors, the {\it vacuum vector} $\1\in V_0$ and
the {\it Virasoro element} $\w\in V_2$.  By definition $Y(\1,z)={\rm
id}_{V}$, and the component operators $\{L(n)|n\in\Z\}$ of
$Y(\w,z)=\sum\limits_{n\in\Z}L(n)z^{-n-2}$ give a representation of
the Virasoro algebra on $V$ with central charge $c.$ Each
homogeneous subspace $V_n\,(n\geq0)$ is an eigenspace for $L(0)$
with eigenvalue $n$.

\begin{definition} A weak $V$-module is a vector space $M$ equipped
with a linear map
$$
\begin{array}{ll}
Y_M: & V \rightarrow ({\rm End}M)[[z,z^{-1}]]\\
 & v \mapsto Y_M(v,z)=\sum\limits_{n \in \Z}v_n z^{-n-1},\ \ v_n \in {\rm End}M
\end{array}
$$
satisfying the following conditions:

1) $v_nw=0$ for $n>>0$ where $v \in V$ and $w \in M$

2) $Y_M( {\textbf 1},z)={\rm Id}_M$

3) The Jacobi identity holds:
\begin{eqnarray}
& &z_0^{-1}\de \left({z_1 - z_2 \over
z_0}\right)Y_M(u,z_1)Y_M(v,z_2)-
z_0^{-1} \de \left({z_2- z_1 \over -z_0}\right)Y_M(v,z_2)Y_M(u,z_1) \nonumber \\
& &\ \ \ \ \ \ \ \ \ \ =z_2^{-1} \de \left({z_1- z_0 \over
z_2}\right)Y_M(Y(u,z_0)v,z_2).
\end{eqnarray}
\end{definition}


\begin{definition}
An admissible $V$-module is a weak $V$-module  which carries a
$\Z_+$-grading $M=\mathop{\bigoplus}\limits_{n \in \Z_+} M(n)$, such
that if $v \in V_r$, then $v_m M(n) \subseteq M(n+r-m-1).$
\end{definition}

\begin{definition}
An ordinary $V$-module is a weak $V$-module which carries a
$\C$-grading $M=\mathop{\bigoplus}\limits_{\la \in \C} M_{\la}$,
such that:

1) $dim(M_{\la})< \infty,$

2) $M_{\la+n}=0$ for fixed $\la$ and $n<<0,$

3) $L(0)w=\la w=\wt(w) w$ for $w \in M_{\la},$ where $L(0)$ is the
component operator of
$Y_M(\omega,z)=\sum\limits_{n\in\Z}L(n)z^{-n-2}.$
\end{definition}

\begin{rem} \ It is easy to see that an ordinary $V$-module is an admissible one. If $W$  is an
ordinary $V$-module, we simply call $W$ a $V$-module.
\end{rem}

We now introduce the Nappi-Witten algebra and some  known related
results.

 The {\it Nappi-Witten Lie algebra} $H_{4}$ is  a $4$-dimensional
vector space
$$
H_{4}=\bC a\oplus \bC b\oplus\bC c\oplus\bC d
$$
equipped with the bracket relations
$$
[a,b]=c,\ \ [d,a]=a,\ \ [d,b]=-b,\ \ [c,d]=[c,a]=[c,b]=0.
$$
and a non-degenerate symmetric bilinear form $(,)$ on $H_{4}$
defined by
$$
(a,b)=1,\quad (c,d)=1,\quad\text{otherwise}, \ (\ ,\ )=0.
$$
The Nappi-Witten algebra $H_4$ is has the following  triangular
decomposition:
$$
H_4=H_4^+\oplus H_4^0\oplus H_4^-=H_4^{\geq 0}\oplus H_4^{-},
$$
where
$$
H_4^+=\bC a,\quad H_4^0=\bC c\oplus\bC d,\quad H_4^-=\bC b.
$$
For $\la\in (H_4^0)^*$,  the  highest weight module (Verma module)
of  $H_4$ is defined by
\begin{equation}
M(\la)=U(H_4)\ot_{U(H_4^{\geq0})}\bC_{\la},\label{Verma}
\end{equation}
where $\bC_{\la}$ is the $1$-dimensional $H_4^{\geq0}$-module, on
which $h\in H_4^0$ acts as multiplication by $\la(h)$ and $H_4^+$
acts as $0$. For convenience, we denote $M(\la)$ by $M(\cc,\d)$ when
$\la(c)=\cc$ and $\la(d)=\d$.

We have the following  lemma from \cite{MP}.

\begin{lem}
For $\cc,\d\in\bC$, $M(\cc,\d)$ is irreducible if and only if
$\cc\neq0$. If $\cc=0$, then the irreducible quotient
$L(\d)=L(0,\d)=\bC v_{\d}$ such that
\begin{equation}
av_{\d}=bv_\d=cv_{\d}=0,\quad dv_{\d}=\d v_{\d}.\label{Trivial}
\end{equation}
\end{lem}
\qed

Let $M=V(\bla,\bm,\bn)$ for $\bla,\bm,\bn\in\bC$, where
$V(\bla,\bm,\bn)$ is the $H_4$-module defined as follows:
 $V(\bla,\bm,\bn)=\mathop{\bigoplus}\limits_{n\in\bZ}
\mathbb{C}v_{n}$ such that
\begin{equation}
dv_{n}=(\bla+n)v_{n},\ cv_{n}=\bm v_{n},\ av_{n}=-\bm v_{n+1},\
bv_{n}=(\bla+\bn+n)v_{n-1}\label{Im}.
\end{equation}
The following result was obtained by M. Willard \cite{Wi}.
\begin{theorem}\label{tt1}
 All irreducible weight modules of
$H_{4}$ with  finite-dimensional weight spaces can be classified
into the following classes:
\begin{itemize}
\item[(1)] Irreducible modules $L(\d)$ for $\d\in \bC$.
\item[(2)] Irreducible highest weight modules $M(\la)$ or irreducible lowest weight
modules $M^-(\mu)$.
\item[(3)] Intermediate series modules $V(\bla,\bm,\bn)$ defined as
(\ref{Im}) such that $\bm\neq 0$ and $\bla+\bn\notin \Z$.

\end{itemize}
\end{theorem}
\qed

 To the pair $(H_{4},(\ ,\ ))$, let $\widehat{H}_{4}$ be the
associated {\it affine Nappi-Witten Lie algebra}  with the
underlying vector space
\begin{equation}
\widehat{H}_{4}=H_{4}\ot\bC[t,t^{-1}]\oplus\bC{\bf k}
\end{equation}
equipped with the bracket relations
\begin{eqnarray}
[h_1\ot t^m,h_2\ot t^n]=[h_1,h_2]\ot
t^{m+n}+m(h_1,h_2)\delta_{m+n,0}{\bf k}, \ [\widehat{H}_{4}, {\bf
k}]=0, \label{DefinitonANW}
\end{eqnarray}
for $h_1,h_2\in H_{4}$ and $m,n\in\bZ$. It is clear that $\H$ has
the $\bZ$-grading:
$$
\H=\coprod_{n\in\bZ}\H^{(n)},
$$
where
$$
\H^{(0)}=H_{4}\oplus\bC\k,\quad\text{and}\quad\H^{(n)}=H_{4}\ot
t^{n}\quad\text{for}\ n\neq0.
$$
Then one has the following graded subalgebras of $\H$:
\begin{eqnarray*}
\H^{(\pm)}=\coprod_{n>0}\H^{(\pm n)},\quad
\H^{(\geq0)}=\coprod_{n\geq0}\H^{(n)}= \H^{(+)}\oplus H_{4}\oplus\bC
\k.
\end{eqnarray*}

Let $M$ be an $H_{4}$-module and $\ell\in{\mathbb C}$. Regard $M$ as
an $\H^{(\geq0)}$-module with $\H^{(+)}$ acting trivially and  $\k$
as the scalar $\ell$. The induced $\H$-module $V_{\H}(\ell, M)$ is
defined by
\begin{equation}
V_{\H}(\ell,M)=\Ind_{\H^{(\geq0)}}^{\H}(M)=U(\H)\ot_{U(\H^{(\geq0)})}M.\label{vermaH4}
\end{equation}

Let $M=L(\d)=\bC v_{\d}$ for $\d\in\bC$ be the one-dimensional
irreducible $H_4$-module defined as in (\ref{Trivial}). Denote
\begin{equation}
V_{\H}(\ell,\mathbbm{d})=U(\H)\ot_{U(\H^{(\geq0)})}\bC
v_{\d}.\label{Vac}
\end{equation}
We have the following result:
\begin{theorem}[\cite{BJP}]\label{Theorem1}
For $\ell,\mathbbm{d}\in\bC$, the $\H$-module
$V_{\H}(\ell,\mathbbm{d})$ is  irreducible  if and only if
$\ell\neq0$. Furthermore, if $\ell=0$, then the irreducible quotient
of $V_{\H}(\ell,\mathbbm{d})$ is isomorphic to the one-dimensional
$\H$-module $\C v_{\d}$.
\end{theorem}

For $h\in H_{4}$, we define the generating function
$$
h(x)=\sum_{n\in\bZ}(h\ot t^n)x^{-n-1}\in\H[[x,x^{-1}]].
$$
Then the defining relations (\ref{DefinitonANW})  can be
equivalently written as
\begin{eqnarray}
[h_1(x_1),h_2(x_2)]=[h_1,h_2](x_2)x_2^{-1}\de\left(\frac{x_1}{x_2}\right)
-(h_1,h_2)\frac{\pa}{\pa
x_1}x_2^{-1}\de\left(\frac{x_1}{x_2}\right)\k.\label{CM}
\end{eqnarray}
Given an $\H$-module $W$, let $h(n)$ denote the operator on $W$
corresponding to $h\ot t^{n}$ for $h\in H_{4}$ and $n\in\bZ$.  We
shall use the notation for  the action of $h(x)$ on $W$:
\begin{equation}
h_W(x)=\sum_{n\in\bZ}h(n)x^{-n-1}\in(\End W)[[x,x^{-1}]].
\end{equation}
\begin{definition}
Let $W$ be a {\it restricted}  $\H$-module  in the sense that for
every  $h\in H_{4}$ and $w\in W$, $h(n)w=0$ for $n$ sufficiently
large. We say that the $\H$-module $W$ is of level $\ell$ if the
central element $\k$ acts as a scalar $\ell$ in $\bC$.
\end{definition}

Let  $\ell$ be a complex number. Consider the induced module defined
as (\ref{Vac})(let $\d=0$):
$$
V_{\H}(\ell,0)=U(\H)\ot_{U(\H^{(\leq 0)})}v_{0}.
$$
 Set
$$
\1=v_0\in V_{\H}(\ell,0).
$$
Then
$$
V_{\H}(\ell,0)=\coprod_{n\geq0}V_{\H}(\ell,0)_{(n)},
$$
where $V_{\H}(\ell,0)_{(n)}$ is spanned by  the vectors
$$
h^{(1)}_{-m_1}\cdots h^{(r)}_{-m_r}\1
$$
for $r\geq0,\ h^{(i)}\in H_{4},\ m_i\geq1$, with $n=m_1+\cdots+m_r$.
It is clear that $V_{\H}(\ell,0)$ is a restricted $\H$-module of
level $\ell$. We can regard $H_{4}$ as a subspace
 of $V_{\H}(\ell,0)$ through the map
 $$
H_{4}\to V_{\H}(\ell,0), \quad h\mapsto h(-1)\1.
 $$
In fact, $H_{4}=V_{\H}(\ell,0)_{(1)}.$

\begin{theorem}[cf. \cite{LL,Lian}]\label{VA} Let $\ell$ be any complex number. Then
there exists  a unique vertex  algebra structure $(V_{\H}(\ell,0),
Y, \1)$ on $V_{\H}(\ell,0)$ such that $ \1 $ is the vacuum vector
and
$$
Y(h,x)=h(x)\in (\End V_{\H}(\ell,0)  )[[x,x^{-1}]]
$$
for $h\in H_{4}$. For $r\geq0, h^{(i)}\in H_{4}, n_i\in\bZ_{+}$, the
vertex operator map for this vertex  algebra structure is given by
\begin{eqnarray*}
Y(h^{(1)}(-n_1)\cdots h^{(r)}(-n_r){\bf1},x)
&=&\NO\pa^{(n_1-1)}h^{(1)}(x)\cdots \pa^{(n_r-1)}h^{(r)}(x)\NO1,
\end{eqnarray*}
where
$$
\partial^{(n)} = \frac{1}{n!}
\left( \frac{d}{dx} \right)^n
$$
$\NO\NO$ is the normal-ordering operation, and $1$ is the identity
operator on $V_{\H}(\ell,0)$.
\end{theorem}

\begin{prop}[cf. \cite{LL}]\label{VAM}
Any module $W$ for the vertex  algebra $V_{\H}(\ell,0)$ is naturally
a restricted  $\H$-module of level $\ell$, with $h_{W}(x)=Y_W(h,x)$
for $h\in H_{4}$. Conversely, any restricted $\H$-module $W$ of
level $\ell$  is naturally a $V_{\H}(\ell,0)$-module as vertex
algebra with
$$
Y_W(h^{(1)}(-n_1)\cdots
h^{(r)}(-n_r){\bf1},x)=\NO\pa^{(n_1-1)}h^{(1)}_W(x)\cdots
\pa^{(n_r-1)}h^{(r)}_W(x)\NO1_W,
$$
for $r\geq0, h^{(i)}\in\H, n_i\in\bZ_{+}$. Furthermore, for any
$V_{\H}(\ell,0)$-module $W$, the $V_{\H}(\ell,0)$-submodules of $W$
coincide with the $\H$-submodules of $W$.
\end{prop}

Let $\ell$ be a non-zero complex number. Set
\begin{eqnarray}
\om&=&\frac{1}{\ell}\left(a(-1)b(-1)\1+c(-1)d(-1)\1\right)-\frac{1}{2\ell}c(-2)\1-\frac{1}{2\ell^2}c(-1)c(-1)\1\label{CV}
\end{eqnarray}
and define operators $L(n)$ for $n\in\bZ$ by
$$
Y(\om,x)=\sum_{n\in\bZ}\om_nx^{-n-1}=\sum_{n\in\bZ}L(n)x^{-n-2}.
$$

Then we have the following results from \cite{BJP}.
\begin{theorem}\label{t2.7}
Let $\ell$ be a complex number such that $\ell\neq0$. Then the
vertex algebra $V_{\H}(\ell,0)$ constructed in Theorem \ref{VA} is a
vertex operator algebra of central charge $4$ with $\om$ defined in
(\ref{CV}) as the conformal vector. The $\bZ$-grading on
$V_{\H}(\ell,0)$ is given by $L(0)$-eigenvalues. Moreover,
$H_{4}=V_{\H}(\ell,0)_{(1)}$, which generates $V_{\H}(\ell,0)$ as a
vertex algebra, and
$$
[h(n),L(m)]=nh(m+n),\quad\text{for}\ h\in H_{4},\ m,n\in\bZ.
$$
\end{theorem}
\begin{theorem}\label{t2.8}
For $\ell\neq0$ and $\d\in\bC$,  the $\H$-module $V_{\H}(\ell,\d)$
is naturally an irreducible ordinary module  for the vertex operator
algebra $V_{\H}(\ell,0)$. Furthermore, the modules $V_{\H}(\ell,\d)$
exhaust all the irreducible ordinary $V_{\H}(\ell,0)$-modules
 up to equivalence.
\end{theorem}
Note that the central charge of the vertex operator algebra
$V_{\H}(\ell,0)$ for $\ell\neq 0$ is independent of $\ell$. We have
the following result.
\begin{theorem}\label{t2.9}
Let  $\ell\neq0$ be any nonzero complex number and let $\sqrt{\ell}$
be one of the two square  roots of $\ell$. Then the vertex operator
algebra $V_{\H}(\ell,0)$ is isomorphic to $V_{\H}(1,0)$ under the
following linear map
$$
\rho: \ V_{\H}(\ell,0)\rightarrow V_{\H}(1,0)
$$
\begin{eqnarray*}
& & d(-m_{41})\cdots d(-m_{4i_{4}})c(-m_{31})\cdots
c(-m_{3i_{3}}) \\
& & a(-m_{21})\cdots a(-m_{2i_{2}})b(-m_{11})\cdots
b(-m_{1i_{1}})\1\\
& \rightarrow &
(\sqrt{\ell})^{i_{1}+i_{2}}\ell^{i_{3}}d(-m_{41})\cdots
d(-m_{4i_{4}})c(-m_{31})\cdots c(-m_{3i_{3}})\\
& & a(-m_{21})\cdots a(-m_{2i_{2}})b(-m_{11})\cdots
b(-m_{1i_{1}})\1,
\end{eqnarray*}
where $i_{k}\geq 0$, $m_{kj}\geq 1$, $m_{k1}\leq m_{k2}\leq \cdots
\leq m_{ki_{k}}$, $k=1,\cdots,4, j=1,2,\cdots,i_{k}$.
\end{theorem}
Following Theorem \ref{t2.9}, for the rest of the paper, we always
assume that $\ell=1$.

\section{Restricted modules for the affine Lie algebra $\widehat{H}_{4}$}

In this section, we will discuss the structures of restricted
$\widehat{H}_{4}$-modules under some conditions which are natural
from the point of view of representations of vertex algebras. We
first have the following lemma.
\begin{lem}\label{l3.0}
Let $W$ be an irreducible $\widehat{H}_{4}$-module, then $c(0)$ and
${\bf k}$ act as scalars on $W$.
\end{lem}
\pf \ Note that $\widehat{H}_{4}$ has countable basis and $W$ is
irreducible. So $W$ is of countable dimension over $\bC$. Then the
lemma follows from the fact that $c(0)$ and ${\bf k}$ are in the
center of $\widehat{H}_{4}$ and ${\rm
Hom}_{\widehat{H}_{4}}(W,W)=\bC$. \qed

\vskip 0.3cm

 Let ${\frak h}=\C c\oplus \C d$ and $\hat{\frak h}={\frak
h}\otimes \C[t,t^{-1}]\oplus \C {\bf k}$ be the associated
subalgebra of the affine Lie algebra $\H$.

 We now introduce a subalgebra
$\widehat{H}_{4}(m,\epsilon)$ of $\widehat{H}_{4}$ for given
$m\in\bZ$ and $\epsilon\in\{0,1\}$ as follows:
\begin{equation}\label{e3.1c}\widehat{H}_{4}(m,\epsilon)=(\bigoplus_{n\geq -m}\bC a\otimes
t^{n})\bigoplus(\bigoplus_{n\geq m+\epsilon}\bC b\otimes
t^{n})\bigoplus((\bC c\bigoplus\bC d)\otimes\bC[t])\bigoplus\bC {\bf
k}. \end{equation} Let $0\neq\al=\d c+\cc d\in{\frak h}^{*}$. We
define a one-dimensional $\widehat{H}_{4}(m,\epsilon)$-module
$\C_{\al}=\C \1_{\al}$ by
$$(a\otimes t^{s})\1_{\al}=(b\otimes t^{r})\1_{\al}=(c\otimes t^{l})\1_{\al}=(d\otimes
t^{l})\1_{\al}=0, \ s\geq -m, \ r\geq m+\epsilon, \ l\geq 1,$$$$
c\1_{\al}=\cc \1_{\al}, \  d{\1}_{\al}=\d \1_{\al}, \ {\bf k}
{\1}_{\al}=\1_{\al}.$$
 Then we get a generalized Verma type
module $V_{\widehat{H}_{4}}(1,(m,\epsilon), \1_{\al})$ as follows:
\begin{equation}\label{eq3.1c}V_{\widehat{H}_{4}}(1,(m,\epsilon),
\1_{\al})=U(\widehat{H}_{4})\otimes_{U(\widehat{H}_{4}(m,\epsilon))}\C_{\al}.
\end{equation} If $m=\cc\in\Z$ and $\epsilon=0$, we briefly denote
$V_{\widehat{H}_{4}}(1,(m,\epsilon), \1_{\al})$ by
$V_{\widehat{H}_{4}}(1,\cc, \1_{\al}).$ Note that
$$V_{\widehat{H}_{4}}(1,0,
\1_{\al})=V_{\widehat{H}_{4}}(1,\d).$$

We first consider the module $V_{\widehat{H}_{4}}(1,\cc, \1_{\al})$
for $\cc\in\Z$. We give the first main result of this section as
follows.
\begin{theorem}\label{t3.1}
Let $0\neq \al=\d c+\cc d\in {\frak h}^{*}$ such that $\cc\in\Z$.
Then the generalized Verma type $\widehat{H}_{4}$-module
$V_{\widehat{H}_{4}}(1,\cc, \1_{\al})$ is irreducible.
\end{theorem}
\pf It is clear that $\sigma:
\widehat{H}_{4}\rightarrow\widehat{H}_{4}$ defined by
$\sigma(a\otimes t^{n})=b\otimes t^{n}$, $\sigma(b\otimes
t^{n})=a\otimes t^{n}$, $\sigma(c\otimes t^{n})=-c\otimes t^{n}$,
$\sigma(d\otimes t^{n})=-d\otimes t^{n}$, $n\in\Z$, $\sigma({\bf
k})={\bf k}$ is an automorphism of $\widehat{H}_{4}$. So we may
assume that $\cc\geq 0$. By PBW theorem and the definition of
$V_{\H}(1,\cc,{\1}_{\al})$,
$$\1_{\al},  d(m_{41})\cdots d(m_{4i_{4}})c(m_{31})\cdots
c(m_{3i_{3}})a(m_{21})\cdots a(m_{2i_{2}})b(m_{11})\cdots
b(m_{1i_{1}})\1_{\al} ,$$ such that $i_{k}\geq 0$, $m_{k1}\leq
m_{k2}\leq \cdots \leq m_{ki_{k}}$, $m_{4j},m_{3j}\leq -1$,
$m_{2j}\leq -\cc-1$, $m_{1j}\leq \cc-1$, $k=1,\cdots,4,
j=1,2,\cdots,i_{k}$, is a basis of $V_{\H}(1,\cc,{\1}_{\al})$.

For $n\in{\mathbb Z}_{+}$, let
$$p_{n}=c(n), \ \ q_{n}=d(-n).$$
Then
$$[p_{m},q_{n}]=m\delta_{m,n}\k.$$
Therefore $${\frak s}=\mathop{\oplus}\limits_{n\in{\mathbb
Z}_{+}}(\C p_{n}\oplus\C q_{n})\oplus\C\k$$ is a Heisenberg algebra
and $V_{\H}(1,\cc,{\1}_{\al})$ is an ${\frak s}$-module such that
$\k$ acts as $1$. Since every highest weight ${\frak s}$-module
generated by one element with $\k$ acting as a non-zero scalar is
irreducible, it follows  that $V_{\H}(1,\cc,{\1}_{\al})$ can be
decomposed into a direct sum of irreducible highest weight modules
of ${\frak s}$ with the highest weight vectors
$$
\1_{\al},  c(m_{31})\cdots c(m_{3i_{3}})a(m_{21})\cdots
a(m_{2i_{2}})b(m_{11})\cdots b(m_{1i_{1}})\1_{\al} ,
$$
such that $i_{k}\geq 0$,  $m_{k1}\leq m_{k2}\leq \cdots \leq
m_{ki_{k}}$, $m_{3j}\leq -1$, $m_{2j}\leq -\cc-1$, $m_{1j}\leq
\cc-1$, $k=1,2,3, j=1,2,\cdots,i_{k}$. For convenience, we denote
the set of such highest weight vectors by ${\cal N}$. Set
$$
\widehat{H}_{4}(\cc,+)=(\bigoplus_{n\geq -\cc}\bC a\otimes
t^{n})\bigoplus(\bigoplus_{n\geq \cc}\bC b\otimes
t^{n})\bigoplus((\bC c\bigoplus\bC d)\otimes\bC t[t])
$$
Let $U$ be a non-zero submodule of $V_{\H}(1,\cc,{\1}_{\al})$. It
suffices to prove that $\1_{\al}\in U$. Let $u$ be a non-zero
element in $U$. Since as an ${\frak s}$-module, $U$ is a direct sum
of irreducible highest weight ${\frak s}$-modules,  we may assume
that
$$
u=\sum_{i=1}^{s}a_{i}c(m_{31}^{(i)})\cdots
c(m_{3i_{3}}^{(i)})a(m_{21}^{(i)})\cdots
a(m_{2i_{2}}^{(i)})b(m_{11}^{(i)})\cdots
b(m_{1i_{1}}^{(i)})b(n_{i1})\cdots b(n_{ir_{i}})\1_{\al},
$$
where $0\neq a_{i}\in\C$, $i_{k}\geq 0$,  $m_{k1}^{(i)}\leq
m_{k2}^{(i)}\leq \cdots \leq m_{ki_{k}}^{(i)}$, $m_{3j}^{(i)}\leq
-1$, $m_{2j}^{(i)}\leq -\cc-1$, $m_{1j}^{(i)}\leq -1$, $0\leq
n_{i1},\cdots,n_{ir_{i}}\leq \cc-1$, $i=1,2,\cdots,s$, $k=1,2,3,
j=1,2,\cdots,i_{k}$, and
$\sum_{j=1}^{1_{1}}m_{1j}^{(1)}\leq\sum_{j=1}^{2_{1}}m_{1j}^{(2)}\leq
\cdots\leq \sum_{j=1}^{s_{1}}m_{1j}^{(s)}$. Let
$$
u^{(i)}=c(m_{31}^{(i)})\cdots
c(m_{3i_{3}}^{(i)})a(m_{21}^{(i)})\cdots
a(m_{2i_{2}}^{(i)})b(m_{11}^{(i)})\cdots
b(m_{1i_{1}}^{(i)})b(n_{i1})\cdots b(n_{ir_{i}})\1_{\al}, \
i=1,2,\cdots,s.$$ Note that
\begin{eqnarray*}
& & a(-m_{11_{1}}^{(1)})\cdots a(-m_{12}^{(1)})
a(-m_{11}^{(1)})u^{(1)}\\
& & =q\prod_{j=1}^{1_{1}}(\cc-m_{1j}^{(1)})c(m_{31}^{(1)})\cdots
c(m_{31_{3}}^{(i)})a(m_{21}^{(1)})\cdots
a(m_{21_{2}}^{(i)})b(n_{i1})\cdots b(n_{ir_{i}})\1_{\al}\\
& & \neq 0,
\end{eqnarray*}
for some $0\neq q\in\Z$ and for the case $b(m_{11}^{(i)})\cdots
b(m_{1i_{1}}^{(i)})\neq b(m_{11}^{(1)})\cdots b(m_{11_{1}}^{(1)})$,
we have
$$
a(-m_{11_{1}}^{(1)})\cdots a(-m_{12}^{(1)})
a(-m_{11}^{(1)})u^{(i)}=0.$$ We may further assume that
$$
u=\sum_{i=1}^{s}a_{i}c(m_{31}^{(i)})\cdots
c(m_{3i_{3}}^{(i)})a(m_{21}^{(i)})\cdots
a(m_{2i_{2}}^{(i)})b(n_{i1})\cdots b(n_{ir_{i}})\1_{\al},
$$
and $$u^{(i)}=c(m_{31}^{(i)})\cdots
c(m_{3i_{3}}^{(i)})a(m_{21}^{(i)})\cdots
a(m_{2i_{2}}^{(i)})b(n_{i1})\cdots b(n_{ir_{i}})\1_{\al}, \
i=1,2,\cdots,s,$$
 where $0\neq a_{i}\in\C$, $i_{k}\geq 0$,
$m_{k1}^{(i)}\leq m_{k2}^{(i)}\leq \cdots \leq m_{ki_{k}}^{(i)}$,
$m_{3j}^{(i)}\leq -1$, $m_{2j}^{(i)}\leq -\cc-1$,  $0\leq n_{i1}\leq
\cdots\leq n_{ir_{i}}\leq \cc-1$, $i=1,2,\cdots,s$, $k=1,2,
j=1,2,\cdots,i_{k}$ and
$\sum_{j=1}^{1_{2}}m_{2j}^{(1)}\leq\sum_{j=1}^{2_{2}}m_{2j}^{(2)}\leq
\cdots\leq \sum_{j=1}^{s_{2}}m_{2j}^{(s)}$. Note also that
\begin{eqnarray*}
& & b(-m_{21_{2}}^{(1)})\cdots b(-m_{22}^{(1)})
b(-m_{21}^{(1)})u^{(1)}\\
& &
=q_{1}\prod_{j=1}^{1_{1}}(-\cc-m_{2j}^{(1)})c(m_{31}^{(1)})\cdots
c(m_{31_{3}}^{(i)})b(n_{i1})\cdots b(n_{ir_{i}})\1_{\al}\\
& & \neq 0,
\end{eqnarray*}
for some $0\neq q_{1}\in \Z$ and for the case $a(m_{21}^{(i)})\cdots
a(m_{2i_{2}}^{(i)})\neq a(m_{21}^{(1)})\cdots b(m_{21_{2}}^{(1)})$,
we have
$$
b(-m_{21_{2}}^{(1)})\cdots b(-m_{22}^{(1)})
b(-m_{21}^{(1)})u^{(i)}=0.$$ We may also assume that
$$
u=\sum_{i=1}^{s}a_{i}c(m_{i1})\cdots c(m_{ip_{i}})b(n_{i1})\cdots
b(n_{ir_{i}})\1_{\al},
$$
where  $0\neq a_{i}\in\C$, $p_{i}, r_{i}\geq 0$, $m_{i1}\leq
m_{i2}\leq \cdots \leq m_{ip_{i}}\leq -1$, $0\leq n_{i1}\leq
\cdots\leq n_{ir_{i}}\leq \cc-1$, $i=1,2,\cdots,s$. Note that
$[d,b]=-b$, so we may assume that
$$r_{1}=r_{2}=\cdots=r_{s}=r.$$
If $r=0$, then by the fact that $[c(p),d(q)]=p\delta_{p+q,0}k$, it
is easy to see that $\1_{\al}\in U$. If $r\neq 0$. Assume that $$
(n_{i1},n_{i2},\cdots,n_{ir})\leq
(n_{i+1,1},n_{i+1,2}\cdots,n_{i+1,r}), \ i=1,2,\cdots,r$$ in the
sense that
$$n_{i1}\leq n_{i2}\leq\cdots\leq n_{ir},$$ and
$$n_{i1}=n_{i+1,1}, \cdots, n_{ik}=n_{i+1,k}, n_{i,k+1}<n_{i+1,k+1},
\ \ {\rm for \ some} \ k.$$
 Then
$$
a(-n_{1r})\cdots a(-n_{11})c(m_{11})\cdots
c(m_{1p_{1}})b(n_{11})\cdots b(n_{1r})=k c(m_{11})\cdots
c(m_{1p_{1}})\1_{\al}\neq 0$$ and for the case that
$(n_{11},\cdots,n_{1r})<(n_{i1},\cdots,n_{ir})$, we have
$$
a(-n_{1r})\cdots a(-n_{11})c(m_{i1})\cdots
c(m_{ip_{i}})b(n_{i1})\cdots b(n_{ir})\1_{\al}=0.
$$
So
$$
a(-n_{1r})\cdots a(-n_{11})u=\sum_{i=1}^{r_{1}}c_{i}c(m_{11})\cdots
c(m_{1p_{1}})\1_{\al}\in U,$$ for some $r_{1}\leq r$ and $c_{i}\neq
0, i=1,2,\cdots,r_{1}$. From this we deduce that $\1_{\al}\in U$.
\qed

\vskip 0.3cm
 Next we consider the second class of  generalized Verma type module  $V_{\H}(1,(m,\epsilon),{\1}_{\al})$ for $\epsilon=1$. Let $0\neq\al=\d c+\cc d\in{\frak h}^{*}$,
$(m,1)\in\Z^{2}$, and let $V_{\H}(1,(m,1),{\1}_{\al})$ be the
generalized Verma type $\H$-module defined by (\ref{eq3.1c}). For
$n\in\Z$, define
\begin{eqnarray*}
&&\deg a(-n-m)=n-1,\ \deg b(-n+m)=n+1,\\
&&\deg c(-n)=n,\  \deg d(-n)=n,\quad\deg\k=0.
\end{eqnarray*}
Then $U(\widehat{H}_4)$ is $\Z$-graded. So we get an $\N$-gradation
of $V_{\H}(1,(m,1),{\1}_{\al})$ by defining $\deg {\1}_{\al}=0$. Let
\begin{equation}\label{e3.2c}
\widehat{H}_{4}((m,1),+)=(\bigoplus_{n\geq -m}\bC a\otimes
t^{n})\bigoplus(\bigoplus_{n\geq m+1}\bC b\otimes
t^{n})\bigoplus((\bC c\bigoplus\bC d)\otimes\bC t[t])
\end{equation}
A homogeneous element $u$ in $V_{\H}(1,(m,1),{\1}_{\al})$ is called
a homogeneous singular vector if
$$
x\cdot u=0,$$  for all $x\in \mathcal {U}(\widehat{H}_{4}((m,1),+)$.
It is obvious that the module generated by a homogeneous singular
vector is a proper submodule. The following is the second main
result of this section.
\begin{theorem}\label{t3.2}
For $0\neq\al=\d c+\cc d\in{\frak h}^{*}$ and $(m,1)\in\Z^{2}$, the
$\H$-module $V_{\H}(1,(m,1),{\1}_{\al})$ is irreducible if and only
if $\cc\notin\Z$. Furthermore,  we have

\begin{itemize}
\item[(i)]  If $\cc=m$, then $u=b(m){\1}_{\al}$ is a
singular vector and $V_{\H}(1,(m,1),{\1}_{\al})/\mathcal {U}(\H)u$
is irreducible and isomorphic to $V_{\widehat{H}_{4}}(1,\cc,
\1_{\al}).$
\item[(ii)] If $\cc+n=m$ for some $n\in\Z_{+}\setminus\{0\}$, then all the
linearly independent singular vectors are
$$
u^{k}=[\sum_{\la\in
P(n)}(a_{\la}c(-\la)b(m)+\sum_{i=1}^{k_{\la}}b_{\la\setminus\la_{i},\la_{i}}c(-\la\setminus\la_{i})b(-\la_{i}+m)]^{k}{\1}_{\al},
$$ for all $k\in\Z_{+}$,  satisfying
\begin{eqnarray*}
&&a_{\la}q_{i}\la_{i}-b_{\la\setminus\la_{i},\la_{i}}=0, \ i=1,2,\cdots,k_{\la},\\
&&b_{\la\setminus\la_{i},\la_{i}}(-n+\la_{i})
+\sum_{j=1}^{k_{\la}}b_{{\la}\setminus\la_{i}\setminus\la_{j},\la_{i}+\la_{j}}=0;
\end{eqnarray*}
\item[(iii)]If $\cc=m+n$ for some $n\in\Z_{+}\setminus\{0\}$, then all the
linearly independent singular vectors are
$$
u^{k}=[\sum_{\la\in
P(n)}\sum_{i=1}^{k_{\la}}c_{\la\setminus\la_{i},\la_{i}}c(-\la\setminus\la_{i})a(-\la_{i}-m)]^{k}{\1}_{\al},$$
for all $k\in\Z_{+}$,  satisfying
$$
q_{i}\la_{i}
c_{\la\setminus\la_{j},\la_{j}}+c_{\la\setminus\la_{i}\setminus\la_{j},\la_{i}+\la_{j}}=0,
\ i,j=1,2,\cdots,k_{\la}, i\neq j,
$$$$\label{dd}c_{\la\setminus\la_{i},\la_{i}}(-n+\la_{i})-\sum_{j=1}^{k_{\la}}c_{\la\setminus\la_{i}\setminus\la_{j},\la_{i}+\la_{j}}=0,
\ i=1,2,\cdots,k_{\la},
$$
where in (ii) and (iii), $P(n)$ is the set of all partitions of
weight $n$,
$$
\la^{(k)}=(\la,\la,\cdots,\la)\in \Z_{+}^{k},\quad
\la=(\la_{1}^{(q_{1})}, \la_{2}^{(q_{2})},\cdots,
\la_{k_{\la}}^{(q_{k_{\la}})})$$ such that
$\sum_{i=1}^{k_{\la}}q_{i}\la_{i}=n$, and
$\la\setminus\la_{i}=(\la_{1}^{(q_{1})},\cdots,
\la_{i}^{(q_{i}-1)},\cdots, \la_{k_{\la}}^{(q_{k_{\la}})})$ and if
$\la=\la_{i}$, then
$c(-\la\setminus\la_{i})b(-\la_{i}+m)=b(-\la_{i}+m)$,
$c(-\la\setminus\la_{i})a(-\la_{i}-m)=a(-\la_{i}-m)$.

\end{itemize}

\end{theorem}

\pf  \ If  $\cc=m$, then it is easy to see that $b(m){\1}_{\al}$
 is a  singular vector. Then (i) follows from Theorem \ref{t3.1}.

We now assume that $\cc\neq m$.

Let ${\frak s}$ be the infinite-dimensional Heisenberg algebra
defined in the proof of Theorem \ref{t3.1}. Then as an ${\frak
s}$-module, $V_{\H}(1,(m,1),{\1}_{\al})$ is a direct sum of
irreducible highest weight ${\frak s}$-modules with $\k$ acting as
$1$ and highest weight vectors in
$$
{\mathcal N}_{\al}=\{b(m)^{i}{\1}_{\al},
c(-\la)a(-\mu)b(-\nu)b(m)^{i}{\1}_{\al},
y_{1}(-\la)y_{2}(-\mu)b(m)^{i}{\1}_{\al},$$$$
x(-\la)b(m)^{i}{\1}_{\al} \ | \ i\geq 0,
(y_{1},y_{2})\in\{(a,b),(c,a),(c,b)\}, x\in\{a,b,c\},
\la,\mu,\nu\in{\mathcal P}\},$$ where  for a partition
$\la=(\la_{1},\cdots,\la_{r})$,
$$d(-\la)=d(-\la_{1})d(-\la_{2})\cdots d(-\la_{r}),$$
$${c}(-\la)=c(-\la_{1})c(-\la_{2})\cdots c(-\la_{r}),$$
$${a}(-\la)=a(-\la_{1}-m)a(-\la_{2}-m)\cdots a(-\la_{r}-m),$$
$${b}(-\la)=b(-\la_{1}+m)b(-\la_{2}+m)\cdots b(-\la_{r}+m).$$

Suppose that $u\in V_{\H}(1,(m,1),{\1}_{\al})$ is a singular vector
such that
$$U(\widehat{H}_4^{+})u=0.$$ Since $V_{\H}(1,(m,1),{\1}_{\al})$ is a direct sum of irreducible highest
weight ${\frak s}$-modules and $u$ is homogeneous, we may assume
that
$$
u=\sum_{i=1}^{3}u_{i},
$$
where
$$u_{1}=\sum_{j=1}^{l_{1}}a_{1j}c(-\la^{(1j)})a(-\mu^{(1j)})b(-\nu^{(1j)})b(m)^{k_{1j}}{\1}_{\al},$$
for
$(\la^{(1j)},\mu^{(1j)},\nu^{(1j)})\succ(\la^{(1,j+1)},\mu^{(1,j+1)},\nu^{(1,j+1)})
, j=1,\cdots, l_{1}-1$,
$$
u_{2}=u_{21}+u_{22}+u_{23},$$ where
$$u_{21}=\sum_{j=1}^{l_{2}}a_{2j}c(-\la^{(2j)})a(-\mu^{(2j)})b(m)^{k_{2j}}{\1}_{\al},
\
u_{22}=\sum_{j=1}^{l_{3}}a_{3j}c(-\la^{(3j)})b(-\nu^{(3j)})b(m)^{k_{3j}}{\1}_{\al},
$$
$$
u_{23}=\sum_{j=1}^{l_{4}}a_{4j}a(-\mu^{(4j)})b(-\nu^{(4j)})b(m)^{k_{4j}}{\1}_{\al},
$$
for $(\la^{(2j)},\mu^{(2j)})\succ(\la^{(2,j+1)},\mu^{(2,j+1)}),
j=1,\cdots, l_{2}-1$;
$(\la^{(3j)},\nu^{(3j)})\succ(\la^{(3,j+1)},\mu^{(3,j+1)}),
j=1,\cdots, l_{3}-1$;
$(\mu^{(4j)},\nu^{(4j)})\succ(\mu^{(4,j+1)},\nu^{(4,j+1)}),
j=1,\cdots, l_{4}-1$, and
$$
u_{3}=u_{31}+u_{32}+u_{33},$$ where
$$
u_{31}=\sum_{j=1}^{l_{5}}a_{5j}c(-\la^{(5j)})b(m)^{k_{5j}}{\1}_{\al},
$$$$
u_{32}=\sum_{j=1}^{l_{6}}a_{6j}a(-\mu^{(6j)})b(m)^{k_{6j}}{\1}_{\al},\
u_{33}=\sum_{j=1}^{l_{7}}a_{7j}b(-\nu^{(7j)})b(m)^{k_{7j}}{\1}_{\al},$$
for $\la^{(5j)}\succ\la^{(5,j+1)}, j=1,\cdots,l_{5}-1$;
$\mu^{(6j)}\succ\mu^{(6,j+1)}, j=1,\cdots,l_{6}-1$;
$\nu^{(7j)}\succ\nu^{(7,j+1)}, j=1,\cdots,l_{7}-1$.

{\bf Case 1} \  $\cc\notin \Z$. One can prove that $u=0$ by the same
method used in the proof of Theorem \ref{t3.1}. Therefore
$V_{\H}(1,(m,1),\al)$ is irreducible.

{\bf Case 2} \ $\cc\in \Z$ and $\cc\neq m$.

Since $au=0$, it is easy to see that $k_{ij}=0$, for $i=4,6,7$. If
$u_{33}\neq 0$, we consider $a(\nu_{1}^{(71)})u$. Since no monomial
of $a(\nu_{1}^{(71)})u_{1}$, $a(\nu_{1}^{(71)})u_{2}$,
$a(\nu_{1}^{(71)})u_{31}$ and $a(\nu_{1}^{(71)})u_{32}$ is of the
form $b(-\eta)$, where $\eta\in{\mathcal P}$, we deduce that
$$
\cc+\nu_{1}^{(71)}-m=0.$$ Similarly, considering
$a(\nu_{r_{7p}}^{(7p)})u$, where
$\nu_{r_{7p}}^{(7p)}=min\{\nu_{r_{7j}}^{(7j)}, \
j=1,2,\cdots,l_{7}\}$, we have
$$
\cc+\nu_{r_{7p}}^{(7p)}-m=0.$$ Then we deduce that $l_{7}=1$ and
$u_{33}=a_{71}(b(-\nu_{1}^{(71)}))^{k_{3}}{\1}_{\al}$ such that
$$\cc+\nu_{1}^{(71)}-m=0, \ \ k_{3}(\nu_{1}^{(71)}+1)=\deg u.$$
If $u_{32}\neq 0$, then by the fact that $b(\mu_{1}^{(61)})u=0$, we
have
$$
\cc-(\mu_{1}^{(61)}+m)=0.$$ Therefore
$$
\mu_{1}^{(61)}+\nu_{1}^{(71)}=0,$$ which is impossible since
$\mu_{1}^{(61)}+\nu_{1}^{(71)}>0$. We deduce that $u_{32}=0$ or
$u_{33}=0$.

{\bf Subcase 1} \ $u_{33}\neq 0$, $u_{32}=0$.

By the fact that
$b(\mu_{1}^{(41)})u=b(\mu_{1}^{(41)})u_{31}=b(\mu_{1}^{(41)})u_{33}=b(\mu_{1}^{(41)})u_{22}=0$,
and both  $b(\mu_{1}^{(41)})u_{1}$ and $b(\mu_{1}^{(41)})u_{21}$
contain no monomials of the forms $a(-\la)b(-\mu)$ and $b(-\mu)$,
where $\la,\mu\in{\mathcal P}$, we have
$$
b(\mu_{1}^{(41)})u_{23}=\sum_{j=1}^{l_{4}}a_{4j}\sum_{p=1}^{r_{4j}}
\delta_{\mu_{1}^{(41)},\mu_{p}^{(4j)}}(-\cc+\mu_{1}^{(41)}+m)$$
$$\cdot a(-\mu^{(4j)})\widehat{a(-\mu_{p}^{(4j)})}b(-\nu^{(4j)}){\1}_{\al}=0,$$
where $\widehat{a(-\mu_{p}^{(4j)})}$ means this factor is deleted.
 So if
$u_{23}\neq 0$, then $-\cc +\mu_{1}^{(41)}+m=0$, which is impossible
since $\cc+\nu_{1}^{(71)}-m=0.$ This proves that $u_{23}=0.$ Let
$Y=b(\mu_{1}^{(21)})$, then $Yu=Yu_{1}+Yu_{21}=0$. If $u_{21}\neq
0$, comparing $Yu_{1}$ and $Yu_{21}$, we have
$$
-\cc +\mu_{1}^{(21)}+m=0,$$ which is not true. So $u_{21}=0$.
Similarly, $u_{1}=0.$ Therefore
$$u=u_{22}+u_{31}+u_{33}.$$
 By the fact that
$a(\nu_{1}^{(31)})u=0$, we can easily deduce that $\nu^{(3j)}\preceq
\nu^{(71)}$, $j=1,\cdots,l_{3}$. Similarly, we have
$\la^{(3j)}\preceq \nu^{(71)},\la^{(5k)}\preceq \nu^{(71)}$,
$j=1,\cdots,l_{3}, k=1,\cdots, l_{5}$. If $\nu_{1}^{(71)}=1$, then
$\cc+1-m=0$, $u_{33}=a_{71}b(-1)^{k_{3}}{\1}_{\al}$ and
$$u_{22}=\sum_{j=1}^{l_{3}}a_{3j}c(-1)^{p_{3j}}b(-1)^{q_{3j}}b(m)^{k_{3j}}{\1}_{\al},\
u_{31}=\sum_{j=1}^{l_{5}}a_{5j}c(-1)^{p_{5j}}b(m)^{k_{5j}}{\1}_{\al}.$$
Since $d(1)u=0$, we have
$$
u_{22}+u_{31}+u_{33}=k(\sum_{j=0}^{k_{3}}C_{k_{3}}^{j}c(-1)^{k_{3}-j}b(-1+m)^{j}b(m)^{k_{3}-j}{\1}_{\al}),$$
for some $k\in \C$. Let $$
u=\sum_{j=0}^{k_{3}}C_{k_{3}}^{j}c(-1)^{k_{3}-j}b(-1+m)^{j}b(m)^{k_{3}-j}{\1}_{\al}=(c(-1)b(m)+b(-1+m))^{k_{3}}{\1}_{\al}.$$
Then we have
$$a(-m)(c(-1)b(m)+b(-1+m))=(c(-1)b(m)+b(-1+m))a(-m).
$$
Thus
$$
a(-m)u=0.$$ It is clear that $a(r-m)u=b(r+m+1)u=c(r+1)u=d(r+1)u=0$,
for $r\geq 0$. We prove that $u$ is a non-zero singular vector for
each $k_{3}\in \N$.

Assume that $\nu_{1}^{(71)}>1$. If $\nu_{1}^{(3j)}=1$ for all
$j=1,2,\cdots,l_{3}$. Then $\nu_{1}^{(71)}=2$. So $\cc+2-m=0$. In
fact, if $\nu_{1}^{(71)}>2$, then
$a(\nu_{1}^{(71)}-1-m)u_{22}=a(\nu_{1}^{(71)}-1-m)u_{31}=0$,
$a(\nu_{1}^{(71)}-1-m)u_{33}\neq 0$, a contradiction. If $k_{3}>1$,
then $a(1-m)u_{31}=0$, and $a(1-m)u_{33}$ contains factor $b(-2+m)$,
but no monomial of $a(1-m)u_{22}$ contains factor $b(-2+m)$. So
$u_{33}=u_{22}=0.$ Then $u_{31}=0$ and we deduce that $u=0$. If
$k_{3}=1$, then $u_{33}=a_{71}b(-2+m){\1}_{\al}$ and one can easily
deduce that
$$
u=( c(-2)b(m)+c(-1)^{2}b(m)+2c(-1)b(-1+m)+2b(-2+m)){\1}_{\al}$$ is a
singular vector. Generally, for $\nu_{1}^{(71)}=2$, we have
$\cc+2-m=0$
 and
$$u=(c(-2)b(m)+c(-1)^{2}b(m)+2c(-1)b(-1+m)+2b(-2+m))^{k}{\1}_{\al}, \ \ k\geq
1.$$

Now assume that $\nu_{1}^{(71)}=n>2$ and $l(\nu^{(71)})=1$. Note
that $\cc+n-m=0$. Let $1\leq p\leq l_{3}$ be  such that
$l(\nu^{(3p)})=max\{l(\nu^{(3j)}), j=1,2,\cdots,l_{3}\} $ and if
$l(\nu^{(3q)})=l(\nu^{(3p)})$, then $\nu^{(3p)}\succeq \nu^{(3q)}$.
If $l(\nu^{(3p)})>1$, then $d(\la^{(3p)})u\neq 0$, a contradiction.
So $l(\nu^{(3j)})=1$ for all $j=1,2,\cdots,l_{3}$. It follows that
$k_{3j}=0, j=1,\cdots,l_{3}$ since $a(-m)u=0$.  For $x\in\Z_{+},
k\in\Z_{+}$, denote $(x,\cdots,x)\in\N^{k}$ by $x^{(k)}$. Then
$$u=\sum_{\la\in
P(n)}[a_{\la}c(-\la)b(m)+\sum_{i=1}^{k_{\la}}b_{\la\setminus\la_{i},\la_{i}}c(-\la\setminus\la_{i})b(-\la_{i})]{\1}_{\al},$$
where $P(n)$ is the set of all partitions of weight $n$,
$\la=(\la_{1}^{(q_{1})}, \la_{2}^{(q_{2})},\cdots,
\la_{k_{\la}}^{(q_{k_{\la}})})$ such that
$\sum_{i=1}^{k_{\la}}q_{i}\la_{i}=n$, and
$\la\setminus\la_{i}=(\la_{1}^{(q_{1})},\cdots,
\la_{i}^{(q_{i}-1)},\cdots, \la_{k_{\la}}^{(q_{k_{\la}})})$ and if
$\la=\la_{i}$, then
$c(-\la\setminus\la_{i})b(-\la_{i}+m)=b(-\la_{i}+m)$. By the fact
that $d(r+1)u=a(r-m)u=b(r+m+1)u=0$ for $r\in\Z_{+}$, and
$\cc+n-m=0$, we deduce  that $a_{\la}$,
$b_{\la\setminus\la_{i},\la_{i}}, i=1,2,\cdots,k_{\la}$ are uniquely
determined by the following equations up to a non-zero scalar.
\begin{equation}
\label{aa} a_{\la}q_{i}\la_{i}-b_{\la\setminus\la_{i},\la_{i}}=0, \
i=1,2,\cdots,k_{\la},
\end{equation}
\begin{equation}\label{bb}
\ b_{\la\setminus\la_{i},\la_{i}}(-n+\la_{i})
+\sum_{j=1}^{k_{\la}}b_{{\la}\setminus\la_{i}\setminus\la_{j},\la_{i}+\la_{j}}=0.
\end{equation}
It is easy to check that the $u$ determined by (\ref{aa}) and
(\ref{bb}) is indeed a non-zero singular vector. Generally
$$u=[\sum_{\la\in
P(n)}(a_{\la}c(-\la)b(m)+\sum_{i=1}^{k_{\la}}b_{\la\setminus\la_{i},\la_{i}}c(-\la\setminus\la_{i})b(-\la_{i}))]^{k}{\1}_{\al},
\ k=1,2,\cdots$$ satisfying (\ref{aa}) and (\ref{bb}) are all the
linearly independent singular vectors.

{\bf Subcase 2} $u_{33}=0$, $u_{32}\neq 0$. Then
$-\cc+\mu_{1}^{(61)}+m=0$.  Similar to the proof for Subcase 1, we
can deduce that $u_{1}=u_{22}=u_{23}=0$, $l_{6}=1$,
$\mu^{(61)}=(\mu_{1}^{(61)},\cdots,\mu_{1}^{(61)})$ and
 $k_{2i}=k_{5j}=0,i=1,2,\cdots,l_{2},
j=1,2,\cdots,l_{5}$.

We first assume that $l(\mu^{(61)})=1$ and $\mu_{1}^{(61)}=n$. Then
it is easy to see that $l(\mu^{(2j)})=1$, $j=1,2,\cdots,l_{2}$,
$u_{21}=0$. Therefore
$$
u=\sum_{\la\in
P(n)}\sum_{i=1}^{k_{\la}}c_{\la\setminus\la_{i},\la_{i}}c(-\la\setminus\la_{i})a(-\la_{i})){\1}_{\al},
$$
where $P(n)$ is the set of all partitions of weight $n$,
$\la^{(k)}=(\la,\la,\cdots,\la)\in \Z_{+}^{k}$,
$\la=(\la_{1}^{(q_{1})}, $ $\la_{2}^{(q_{2})},\cdots,
\la_{k_{\la}}^{(q_{k_{\la}})})$ such that
$\sum_{i=1}^{k_{\la}}q_{i}\la_{i}=n$, and
$\la\setminus\la_{i}=(\la_{1}^{(q_{1})},\cdots,
\la_{i}^{(q_{i}-1)},\cdots, \la_{k_{\la}}^{(q_{k_{\la}})})$. By
$d(n)u=b(n)u=0$ for $n\geq 1$, we deduce that
\begin{equation}
\label{cc} q_{i}\la_{i}
c_{\la\setminus\la_{j},\la_{j}}+c_{\la\setminus\la_{i}\setminus\la_{j},\la_{i}+\la_{j}}=0,
\ i,j=1,2,\cdots,k_{\la}, i\neq j,
\end{equation}
\begin{equation}\label{ee}c_{\la\setminus\la_{i},\la_{i}}(-n+\la_{i})-\sum_{j=1}^{k_{\la}}c_{\la\setminus\la_{i}\setminus\la_{j},\la_{i}+\la_{j}}=0,
\ i=1,2,\cdots,k_{\la}.
\end{equation}
Actually, $u$ is uniquely determined by (\ref{cc}) and (\ref{ee}).

Next for the case that  $l(\mu^{(61)})>1$ and $\mu_{1}^{(61)}=n$,
one can easily deduce that
$$
u=[\sum_{\la\in
P(n)}\sum_{i=1}^{k_{\la}}c_{\la\setminus\la_{i},\la_{i}}c(-\la\setminus\la_{i})a(-\la_{i}))]^{k}{\1}_{\al},\
k\in\Z_{+}$$ satisfying (\ref{cc}) and (\ref{ee}) for some $k\in\N.
$

Conversely, for $k\in\Z$, $k\geq 1$, it is easy to check that
$$u=[\sum_{\la\in
P(n)}\sum_{i=1}^{k_{\la}}c_{\la\setminus\la_{i},\la_{i}}c(-\la\setminus\la_{i})a(-\la_{i}))]^{k}{\1}_{\al},\
k\in\Z_{+}$$ satisfying (\ref{cc}) and (\ref{ee}) is a singular
vector. The proof of the theorem is complete. \qed
\begin{coro}\label{co3.1}
(1) \ If $\cc+n=m$ for some $n\in\Z_{+}$, then
$$L_{\H}(1,(m,1),{\1}_{\al})=V_{\H}(1,(m,1),{\1}_{\al})/U(\H)u^{1}$$ is irreducible, where
$$u^{1}=\sum_{\la\in
P(n)}(a_{\la}c(-\la)b(m)+\sum_{i=1}^{k_{\la}}b_{\la\setminus\la_{i},\la_{i}}c(-\la\setminus\la_{i})b(-\la_{i})){\1}_{\al},
$$ satisfying (\ref{aa}) and (\ref{bb});

(2) If $\cc-n=m$ for some $n\in\Z_{+}$, then
$$L_{\H}(1,(m,1),{\1}_{\al})=V_{\H}(1,(m,1),{\1}_{\al})/U(\H)u^{1}$$ is irreducible, where
$$
u^{1}=\sum_{\la\in
P(n)}\sum_{i=1}^{k_{\la}}c_{\la\setminus\la_{i},\la_{i}}c(-\la\setminus\la_{i})a(-\la_{i})){\1}_{\al},$$
satisfying (\ref{cc}) and (\ref{ee}).
\end{coro}

\pf We only prove (1), since the proof of (2) is quite similar. Let
$$V'={U}(\H)u^{1}.$$
 By Theorem \ref{t3.2} (ii), all the homogeneous singular vectors
$$u^{k}=[\sum_{\la\in
P(n)}(a_{\la}c(-\la)b(m)+\sum_{i=1}^{k_{\la}}b_{\la\setminus\la_{i},\la_{i}}c(-\la\setminus\la_{i})b(-\la_{i}))]^{k}{\1}_{\al},
\ k=1,2,\cdots$$ lie in $V'$. Now let $U/V'$ be a proper submodule
of $V_{\H}(1,(m,1),{\1}_{\al})/V'$, where $U$ is a proper submodule
of $V_{\H}(1,(m,1),{\1}_{\al})$. Recall that
$$V_{\H}(1,(m,1),{\1}_{\al})=\bigoplus_{j=1}^{\infty}V_{\H}(1,(m,1),{\1}_{\al})_{j}.$$ Let
$$
v=\sum\limits_{i}v_{i}\in U, $$ where $0\neq v_{i}\in
V_{\H}(1,(m,1),{\1}_{\al})_{i},$ be such that $|v|=\sum\limits_{i}i$
is the minimal. Then
$$
\widehat{H}_{4}((m,1),+)v_{i}=0,$$ for each $i$, where
$\widehat{H}_{4}((m,1),+)$ is the same as (\ref{e3.2c}). By Theorem
\ref{t3.2} (ii), $v_{i}\in V'$. So $U=V'$. This proves that
$V_{\H}(1,(m,1),{\1}_{\al})/V'$ is irreducible. \qed

Finally we introduce the so-called third class of generalized Verma
type modules of $\H$.  For $m\in\Z$, let $\widehat{H}_{4}(m,0)$ be
the subalgebra of $\widehat{H}_{4}$ defined by (\ref{e3.1c}) for the
case $\epsilon=0$, that is,
$$\widehat{H}_{4}(m,0)=(\bigoplus_{n\geq -m}\bC a\otimes
t^{n})\bigoplus(\bigoplus_{n\geq m}\bC b\otimes t^{n})\bigoplus((\bC
c\bigoplus\bC d)\otimes\bC[t])\bigoplus\bC {\bf k}.$$
 Set
\begin{equation}\label{e3.3}
\widehat{H}_{4}(m,0)^{+}=(\bigoplus_{n\geq -m+1}\bC a\otimes
t^{n})\bigoplus(\bigoplus_{n\geq m+1}\bC b\otimes
t^{n})\bigoplus((\bC c\bigoplus\bC d)\otimes\bC t[t]).
\end{equation} For $\bla, \bm, \bn \in \C$,  let
$$V(\bla,
\bm,\bn,m)=\bigoplus_{n\in\Z}\C v_{n}$$ be an infinite dimensional
module of $\widehat{H}_{4}(m,0)$ defined by
$$
\widehat{H}_{4}(m,0)^{+}v_{n}=0, \ {\bf k}\cdot v_{n}=v_{n} \ \
n\in\Z,
$$
$$
c(0)v_{n}=(\bm+m) v_{n}, \ a(-m)v_{n}=-\bm v_{n+1},
$$
$$
d(0)v_{n}=(\bla+n)v_{n}, \ b(m)v_{n}=(\bla+\bn+n)v_{n-1}.$$ The
following lemma is obvious.
\begin{lem}
$V(\bla, \bm, \bn, m)$ is an irreducible
$\widehat{H}_{4}(m,0)$-module if and only if $\bm\neq 0$ and
$\bla+\bn\notin \Z.$
\end{lem}
Let
 $$ V_{\H}(\bla,\bm,\bn,m)={U}(\H)\otimes_{{U}(\widehat{H}_{4}(m,0))}V(\bla,\bm,\bn,m)$$
be the induced $\H$-module.

For $n\in\Z$, define
\begin{eqnarray*}
&&\deg a(-n-m)=n,\ \deg b(-n+m)=n,\\
&&\deg c(-n)=n,\  \deg d(-n)=n,\quad\deg\k=0.
\end{eqnarray*}
Then $U(\widehat{H}_4)$ is $\Z$-graded. So we get an $\N$-gradation
of $V_{\H}(\bla,\bm,\bn,m)$ by defining $\deg v_{n}=0$ for $n\in\Z$.
A homogeneous vector $u$ in $V_{\H}(\bla,\bm,\bn,m)$ is called a
homogeneous singular vector if
$$
\widehat{H}_{4}(m,0)^{+}(u)=0.
$$

 Similar to Theorem \ref{t3.2}, we have the  following last main result of this section.

\begin{theorem}\label{t3.3}Let $\bla,\bm,\bn\in\bC$ be such that $\bla+\bn\notin \bZ$, $\bm\neq 0$. Then  the $\H$-module
$V_{\H}(\bla,\bm,\bn,m)$ is irreducible if and only if $\bm+n\neq 0$
for all $n\in \mathbb{Z}$. Furthermore, we have
\begin{itemize}
\item[(i)]If $\bm+n=0$ for some $n\in\Z_{+}$, then all the linearly independent homogeneous singular vectors are
$$
u^{(k,j)}=[\sum_{\la\in
P(n)}(a_{\la}c(-\la)b(m)+\sum_{i=1}^{k_{\la}}b_{\la\setminus\la_{i},\la_{i}}c(-\la\setminus\la_{i})b(-\la_{i}+m))]^{k}v_{j},
\ j\in\Z, \ k=1,2,\cdots $$  satisfying
\begin{equation}\label{aa1}a_{\la}q_{i}\la_{i}-b_{\la\setminus\la_{i},\la_{i}}=0, \
i=1,2,\cdots,k_{\la}, \end{equation}
\begin{equation}\label{bb1}
b_{\la\setminus\la_{i},\la_{i}}(\bm+\la_{i})+\sum_{j=1}^{k_{\la}}b_{{\la}\setminus\la_{i}\setminus\la_{j},\la_{i}+\la_{j}}=0;
\end{equation}
The submodules generated by $u^{(k,i)}$ and $u^{(k,j)}$ for
$i,j\in\Z$, $i\neq j$ are isomorphic to each other.

\item[(ii)] If $\bm-n=0$ for some $n\in\Z_{+}$, then
all the homogeneous singular vectors are
$$
u^{(k,j)}=[\sum_{\la\in
P(n)}\sum_{i=1}^{k_{\la}}c_{\la\setminus\la_{i},\la_{i}}c(-\la\setminus\la_{i})a(-\la_{i}))]^{k}v_{j},
\ j\in\Z,  \ k=1,2,\cdots $$
 satisfying
\begin{equation}\label{cc1}
q_{i}\la_{i}
c_{\la\setminus\la_{j},\la_{j}}+c_{\la\setminus\la_{i}\setminus\la_{j},\la_{i}+\la_{j}}=0,
\ i,j=1,2,\cdots,k_{\la}, i\neq j,
\end{equation}
\begin{equation}\label{ee1}c_{\la\setminus\la_{i},\la_{i}}(-\bm+\la_{i})-\sum_{j=1}^{k_{\la}}
c_{\la\setminus\la_{i}\setminus\la_{j},\la_{i}+\la_{j}}=0, \
i=1,2,\cdots,k_{\la},
\end{equation}
and  the submodules generated by $u^{(k,i)}$ and $u^{(k,j)}$ for
$i,j\in\Z$, $i\neq j$ are isomorphic to each other, where in (i) and
(ii), $P(n)$ is the set of all partitions of weight $n$,
$$\la^{(k)}=(\la,\la,\cdots,\la)\in \Z_{+}^{k}, \
\la=(\la_{1}^{(q_{1})},  \ \la_{2}^{(q_{2})},\cdots,
\la_{k_{\la}}^{(q_{k_{\la}})})$$ such that
$\sum_{i=1}^{k_{\la}}q_{i}\la_{i}=n$, and
$\la\setminus\la_{i}=(\la_{1}^{(q_{1})},\cdots,
\la_{i}^{(q_{i}-1)},\cdots, \la_{k_{\la}}^{(q_{k_{\la}})})$ and if
$\la=\la_{i}$, then
$c(-\la\setminus\la_{i})b(-\la_{i})=b(-\la_{i})$,
$c(-\la\setminus\la_{i})a(-\la_{i})=a(-\la_{i})$.

\end{itemize}
\end{theorem}
\begin{coro}\label{co3.2}
(1) \ If $\bm+n=m$ for some $n\in\Z_{+}$, then
$$L_{\H}(\bla,\bm,\bn,m)=V_{\H}(\bla,\bm,\bn,m)/U(\H)u^{1}$$ is irreducible, where
$$u^{1}=\sum_{\la\in
P(n)}(a_{\la}c(-\la)b(m)+\sum_{i=1}^{k_{\la}}b_{\la\setminus\la_{i},\la_{i}}c(-\la\setminus\la_{i})b(-\la_{i}))v_{0},
$$ satisfying (\ref{aa1}) and (\ref{bb1});

(2) If $\bm-n=m$ for some $n\in\Z_{+}$, then
$$L_{\H}(\bla,\bm,\bn,m)=V_{\H}(\bla,\bm,\bn,m)/U(\H)u^{1}$$ is irreducible, where
$$
u^{1}=\sum_{\la\in
P(n)}\sum_{i=1}^{k_{\la}}c_{\la\setminus\la_{i},\la_{i}}c(-\la\setminus\la_{i})a(-\la_{i}))v_{0},$$
satisfying (\ref{cc1}) and (\ref{ee1}).
\end{coro}

\section{Extension of the vertex operator algebra
$V_{\widehat{H}_{4}}(1,0)$}

  Set
$${\frak h}=\bC c(-1)\1 \oplus\bC d(-1)\1.$$ Then
$$L(n)h=\delta_{n,0}h,\quad h(n)h^{'}=(h,h^{'})\delta_{n,1}\1 \quad \text{for}\ n\in\bZ_{+},\quad h,h^{'}\in {\frak h} $$
and $( \ , \ )$ restricted to ${\frak h}$ is a non-degenerate
symmetric bilinear form. Thus we can identify ${\frak h}$ with its
dual ${\frak h}^{\ast}$. One can easily see that for any $h\in
{\frak h}$, $h(0)$ acts semisimply on $V_{\widehat{H}_{4}}(1,0)$. So
$$V_{\widehat{H}_{4}}(1,0)=\mathop{\oplus}\limits_{\alpha\in {\frak h}}V_{\widehat{H}_{4}}(1,0)^{(0,\alpha)},$$
where $$V_{\widehat{H}_{4}}(1,0)^{(0,\alpha)} =\{v\in
V_{\widehat{H}_{4}}(1,0) | h(0)v=(\alpha,h)v \quad \text{for}\ h\in
{\frak h} \}.$$ Set $$P=\{\alpha\in {\frak
h}|V_{\widehat{H}_{4}}(1,0)^{(0,\alpha)} \neq 0\}$$ According to
\cite{Li2}, $P$ equipped with the bilinear form (\ ,\ ) is a
lattice. Then
$$P=\bZ c$$
Set $$L=\bZ c \oplus \bZ d$$
Then for any $\alpha \in L$, $\alpha(0)$
acting on $V_{\widehat{H}_{4}}(1,0)$ has only integral eigenvalues.
This amounts to $L \subset P^{0}$, where
$$P^{0}=\{h\in {\frak h}|(h,\alpha)\in \bZ \quad \text{for}\ \alpha \in P\}$$
Let $\epsilon$: $L\times L \ \rightarrow \{ \pm 1\}$ be a 2-cocycle
on $L$ such that
$$\epsilon(c,c)=1, \quad \epsilon(c,d)=-1, \quad \epsilon(d,d)=1,$$
and
$$\epsilon(\alpha,\beta)\epsilon(\beta,\alpha)^{-1}
=(-1)^{(\alpha,\beta)+(\alpha,\alpha)(\beta,\beta)}=(-1)^{(\alpha,\beta)}
\quad \text{for}\ \alpha,\beta\in L.$$

For $\al\in L$, let
$$E^{\pm}(\alpha,z)=exp\left(\sum^{\infty}_{n=1}\frac{\alpha(\pm n)}
{\pm n}z^{\mp n}\right),\quad
\Delta(\alpha,z)=z^{\alpha(0)}E^{+}(-\alpha,-z).$$ Let $(W, Y(\cdot,
z))$ be any (irreducible) weak $V_{\widehat{H}_{4}}(1,0)$-module,
then by Proposition 2.9 of \cite{Li2},  $(W^{(\al)}, Y_{\al}(\cdot,
z))=(W, Y(\bigtriangleup(\al, z)\cdot, z))$ carries the structure of
an (irreducible) weak $V_{\widehat{H}_{4}}(1,0)$-module. We expand
$Y_{\al}(u, z)$ for $u\in V_{\widehat{H}_{4}}(1,0)$ as follows:
$$
Y_{\al}(u, z))=\sum\limits_{n\in\Z}u_{(\al,n)}z^{-n-1}.$$ We have
the following result. \vskip 0.3cm
\begin{prop}\label{p4.1}
 (1) \ Let $k\in\Z$ and $\al=kc\in L$, then
$V_{\widehat{H}_{4}}(1,0)^{(\al)}\cong V_{\widehat{H}_{4}}(1,\d)$
with $\d=k$ and  $V_{\widehat{H}_{4}}(1,0)^{(\al)}$ is a simple
current of $V_{\widehat{H}_{4}}(1,0)$.

(2) Let $m_{1},m_{2}\in\Z, m_{2}\neq 0$ and $\al=m_{1}c+m_{2}d$.
Then $V_{\widehat{H}_{4}}(1,0)^{(\alpha)}$ is isomorphic to the
generalized Verma type module $V_{\widehat{H}_{4}}(1,m_{2},
\1_{\al}).$

(3)   $V_{\widehat{H}_{4}}(1,0)^{(\al)}$ is irreducible for each
$\al\in L$ and for $\al,\be\in L$,
$V_{\widehat{H}_{4}}(1,0)^{(\al)}\cong
V_{\widehat{H}_{4}}(1,0)^{(\be)}$ if and only if $\al=\be$.
\end{prop}

\pf For $\al=kc$, let $W$ be any irreducible weak
$V_{\widehat{H}_{4}}(1,0)$-module. By Lemma \ref{l3.0}, $c(0)$ acts
as a scalar on $W$. Then by Theorem 2.15 of \cite {Li3} ( see also
Theorem 2.13 of \cite{Li2}), $V_{\widehat{H}_{4}}(1,0)^{(\al)}$ is a
simple current of  $V_{\widehat{H}_{4}}(1,0)$. By Theorem
\ref{t2.8}, $W\cong V_{\widehat{H}_{4}}(1,\d)$ for some $\d\in\bC$.
From the definition of $Y_{\al}(\cdot, z)=Y(\bigtriangleup(\al,
z)\cdot, z)$, it is easy to check that $W^{(\al)}\cong
V_{\widehat{H}_{4}}(1,\d+k)$. Therefore (1) holds.

For (2), by the definition of $Y_{\al}$, we have
$$
Y_{\al}(a(-1){\bf 1},z){\bf 1}=z^{m_{2}}Y(a(-1){\bf 1}, z){\bf 1},
 \ \
Y_{\al}(b(-1){\bf 1},z){\bf 1}=z^{-m_{2}}Y(b(-1){\bf 1}, z){\bf 1},
$$
$$
Y_{\al}(c(-1){\bf 1},z){\bf 1}=Y(c(-1){\bf 1}+z^{-1}m_{2}{\bf 1},
z){\bf 1}, $$$$Y_{\al}(d(-1){\bf 1},z){\bf 1}=Y(d(-1){\bf
1}+z^{-1}m_{1}{\bf 1}, z){\bf 1},
$$
So we have
$$(a(-1){\bf 1})_{(\al,n)}{\bf 1}=0, \ \ (b(-1){\bf 1})_{(\al,m)}{\bf 1}=0, \ {\rm
for} \ \ n\geq -m_{2}, \ m\geq m_{2},$$
$$
(c(-1){\bf 1})_{(\al,n)}{\bf 1}=(d(-1){\bf 1})_{(\al,n)}{\bf 1}=0, \
\ {\rm for} \ \ n\geq 1,$$
$$
(c(-1){\bf 1})_{(\al,0)}{\bf 1}=m_{2}{\bf 1}, \ \ (d(-1){\bf
1})_{(\al,0)}{\bf 1}=m_{1}{\bf 1}.$$ Therefore as an
$\widehat{H}_{4}$-module, $V_{\widehat{H}_{4}}(1,0)^{(\alpha)}$ is a
quotient of the generalized Verma type module
$V_{\widehat{H}_{4}}(1, m_{2}, \1_{\al})$. By Theorem \ref{t3.1},
$V_{\widehat{H}_{4}}(1,0)^{(\alpha)}$ is isomorphic to
$V_{\widehat{H}_{4}}(1, m_{2}, \1_{\al})$.

(3) follows from Theorem \ref{t3.1} and (1)-(2).
 \qed

\vskip 0.3cm Furthermore, we have
\begin{prop}\label{p4.2} Let $\al=m_{1}c+m_{2}d\in L, \
\be=n_{1}c+n_{2}d\in {\frak h}^{*}$ and $m\in \Z$. Then

 (1) $$V_{\H}(1,(m,1),{\1}_{\be})^{(\al)}\cong V_{\H}(1,(m+m_{2},1),{\1}_{\al+\be}), $$
 $$\ L_{\H}(1,(m,1),{\1}_{\be})^{(\al)}\cong L_{\H}(1,(m+m_{2},1),{\1}_{\al+\be}).$$

 (2) Let $\bla,\bm,\bn\in\bC$ be such that $\bla+\bn\notin \bZ$, $\bm\neq
 0$. Then
$$V_{\H}(\bla,\bm,\bn,m)^{(\al)}\cong
V_{\H}(\bla+m_{1},\bm,\bn,m+m_{2}),$$
$$L_{\H}(\bla,\bm,\bn,m)^{(\al)}\cong
L_{\H}(\bla+m_{1},\bm,\bn,m+m_{2}).$$
\end{prop}

 For convenience, we give the following definition.
\begin{definition} Let ${\cal F}$ be the category of irreducible restricted $\H$-modules
$W$ of level 1 such that $d(0)$ is semisimple, and there exists an
integer $m$ such that  $h(n), a(-m+n)$ and $b(m+n)$ for $n\geq 1$,
$h\in {\frak h}$ are locally nilpotent on $W$.
\end{definition}
 \vskip 0.3cm

The following lemma is easy to prove.
\begin{lem}\label{l3.2}
Let $W\in{\cal F}$. Then there exist $m\in\Z$ and  a non-zero
element $w$ in $W$ such that
$$
h(n)w=a(-m+n)w=b(m+n)w=0, \ \forall \ n\geq 1, \ \ \forall \ h\in
{\frak h}.$$
\end{lem}

\vskip 0.3cm

We are now in a position to state  the first  main result of this
section.
\begin{theorem}\label{tt4.2} Let $W\in{\cal F}$ be such that
$c(0)$  acts on $W$ as  $\cc$. Denote $\al=\d c+\cc d$, then $W$ is
isomorphic to one of the following irreducible modules:
$$V_{\H}(1,m,{\1}_{\al}), V_{\H}(1,(m,1),{\1}_{\al}), \  V_{\H}(\bla,\bm,\bn,m)\
L_{\H}(1,(m,1),{\1}_{\al}), \ L_{\H}(\bla,\bm,\bn,m)$$
 given in Theorems \ref{t3.1}, \ref{t3.2} and
\ref{t3.3}, Corollary \ref{co3.1} and Corollary \ref{co3.2}
respectively.
\end{theorem}

\pf Let $W\in {\cal F}$, by the definition of ${\cal F}$ and Lemma
\ref{l3.2}, there exist an integer $m$ and a non-zero element $u$ in
$W$ such that
$$
c(0)u=\cc u, \ d(0)u=\d u,
$$
for some $\d\in \C$, and
$$ h(n)u=a(-m+n)u=b(m+n)u=0, \ \forall \
n\geq 1, \ \ \forall \ h\in {\frak h}.$$

We first assume that $m=0$. Then \begin{equation}\label{ee4}
h(n)u=a(n)u=b(n)u=0, \ \forall \ n\geq 1, \ \ \forall \ h\in {\frak
h}.
\end{equation}
 Let $W^{0}$ be the subspace of $W$ consisting of elements
 satisfying (\ref{ee4}). Then $W^{0}$ is a module of $H_{4}$. Since
 $W$ is irreducible, it follows that $W^{0}$ is an irreducible
 $H_{4}$-module. By Theorem \ref{tt1}, Theorems \ref{t3.1},
 \ref{t3.2} and \ref{t3.3},
 Corollaries
 \ref{co3.1} and  \ref{co3.2}, we have
  $$W\cong V_{\H}(1,\d), \ \ {\rm if} \ \cc=0,$$
 $$W\cong V_{\H}(1,(0,1),{\1}_{\al}) \ {\rm or} \
 V_{\H}(\bla,\bm,\bn,0), \ \
{\rm  if} \  0\neq \cc\notin\Z,$$ and $$W\cong
L_{\H}(1,(0,1),{\1}_{\al}) \ {\rm or} \ L_{\H}(\bla,\bm,\bn,0), \ \
{\rm if} \ 0\neq \cc\in\Z.$$

\vskip 0.3cm
 Next we assume that $m\neq 0$. Note that  $W$ can be viewed as a module for the
vertex algebra $V_{\H}(1,0)$. Let $(\cdot,z)$ be the corresponding
operator.  Let $\beta=-md$. By the proof of Proposition \ref{p4.1},
$(W^{(\beta)}=W,Y_{\beta}(\cdot,z))$ is a module of $V_{\H}(1,0)$
satisfying
$$
(x(-1){\bf 1})_{(\beta,n)}u=0, \ \forall \ x\in H_{4}, \ \forall \
n\in\Z_{+},
$$
$$
(c(-1){\bf 1})_{(\beta,0)}u=(\cc-m)u, \ (d(-1){\bf
1})_{(\beta,0)}u=\d u.$$ Therefore $W^{(\beta)}$ is isomorphic to
one of the following module:
$$V_{\H}(1,0,{\1}_{\al+\beta}), \ V_{\H}(1,(0,1),{\1}_{\al+\beta}), \  V_{\H}(\bla,\bm,\bn,0), \
L_{\H}(1,(0,1),{\1}_{\al+\beta}), \ L_{\H}(\bla,\bm,\bn,0).$$ This
means that $W$ is isomorphic to one of the irreducible modules:
$$V_{\H}(1,m,{\1}_{\al}), \ V_{\H}(1,(m,1),{\1}_{\al}), \  V_{\H}(\bla,\bm,\bn,m), \
L_{\H}(1,(m,1),{\1}_{\al}), \ L_{\H}(\bla,\bm,\bn,m).$$

 \qed

Set
$$V_{\widehat{H}_{4}}(1,0)[L]=\bC[L]\otimes
V_{\widehat{H}_{4}}(1,0).$$  For $\alpha,\beta \in L, \ u,v \in
V_{\widehat{H}_{4}}(1,0),$ we define a vertex operator map $Y$ on
$V_{\widehat{H}_{4}}(1,0)[L]$ by
$$Y(e^{\alpha}\otimes u,z)(e^{\beta}\otimes
v)=\epsilon(\alpha,\beta)e^{\alpha+\beta}\otimes z^{
(\alpha,\beta)}E^{-}(-\alpha,z)Y(\Delta(\beta,z)u,z)E^{+}(-\alpha,z)(-z)^{\alpha(0)}v.$$

From \cite{Li2}, we have

\begin{theorem} $V_{\widehat{H}_{4}}(1,0)[L]$ equipped with the vertex
operator map $Y$ is a  vertex algebra.\end{theorem}

\vskip 0.3cm

 It is easy to see that for each $\al\in L$, $\bC
e^{\al}\otimes V_{\widehat{H}_{4}}(1,0)$ is isomorphic to
$V_{\widehat{H}_{4}}(1,0)^{(\al)}$ as
$V_{\widehat{H}_{4}}(1,0)$-modules. So we also denote
$$V_{\widehat{H}_{4}}(1,0)^{(\al)}=\bC e^{\al}\otimes V_{\widehat{H}_{4}}(1,0).$$

 It is obvious that  $L(0)$ has
integral eigenvalues on $V_{\widehat{H}_{4}}(1,0)[L]$. In fact, for
$$v=e^{k_{1}c+k_{2}d}\otimes d(-m_{11})\cdots
d(-m_{1i_{1}})c(-m_{21})\cdots c(-m_{2i_{2}}) a(-m_{31})$$$$\cdots
a(-m_{3i_{3}})b(-m_{41})\cdots b(-m_{4i_{4}}){\bf 1},$$
$$L(0)v=(k_{1}k_{2}+\sum_{k=1}^{4}\sum_{j=1}^{i_{k}}m_{kj}+k_{2}(i_{3}-i_{4}))v.$$
 Moreover, the same proof of Theorem 3.5 in \cite{DLM1} gives that for $v\in
V_{\widehat{H}_{4}}(1,0)[L]$,
$$[L(-1), Y(v,z)]=\frac{d}{dz}Y(v,z).$$Thus
$V_{\widehat{H}_{4}}(1,0)[L]$ is a vertex operator algebra except
the two grading restrictions. \vskip 0.3cm

\begin{theorem} Let $W$ be a module for the vertex  algebra
$V_{\widehat{H}_{4}}(1,0)$ such that for each $\al\in L$, $\al(0)$
acts on $W$ semisimply with only integral eigenvalues.  Set
$$W[L]=\bC [L]\otimes W.$$ For $\al,\be\in L$, $v\in
V_{\widehat{H}_{4}}(1,0)$, $w\in W$, we define
$$
Y_{W[L]}(e^{\al}\otimes v,z)(e^{\be}\otimes w)=\epsilon
(\al,\be)e^{\al+\be}\otimes
z^{(\al,\be)}E^{-}(-\al,z)Y_{W}(\bigtriangleup(\be,z)v,z)E^{+}(-\al,z)(-z)^{\al(0)}w.$$
Then $W[L]$ carries the structure of a
$V_{\widehat{H}_{4}}(1,0)[L]$-module. Furthermore, if $W$ is
irreducible, then $W[L]$ is irreducible. In particular,
$V_{\widehat{H}_{4}}(1,0)[L]$ is a simple vertex operator algebra.

\end{theorem}
\begin{lem}\label{l4.3}
 Let $U$ be an irreducible
$V_{\widehat{H}_{4}}(1,0)[L]$-module. Then

(1) There exists an element $w$ in $U$ such that for all $\al\in L$
and $n\geq 1$,
$$\al(n)w=0.$$

(2) For $\al\in L$, $\al(0)$ acts on $U$ semisimply with only
integral eigenvalues and all $\al(n)$ for $n\geq 1$ are locally
nilpotent.
\end{lem}

\pf Note that $V_{\widehat{H}_{4}}(1,0)[L]$ contains the lattice
vertex algebra $V_{L}$ and $V_{\widehat{H}_{4}}(1,0)[L]$ is simple.
By Lemma 3.15 in \cite{DLM2}, there exists an element $w$ in $U$
such that for all $\al\in L$, and $n\geq 1$,
$$\al(n)w=0.$$
  Let $\al$ be any element in $L$. Since $V_{\widehat{H}_{4}}(1,0)[L]$
is simple, similar to the argument in \cite{D} (see also \cite{LW}),
from \cite{DL}, we have $Y(e^{\al}, z)w\neq 0$. Then there exists
$k\in\Z$ such that $e^{\al}(k)w\neq 0$ and $e^{\al}(n)w=0$ for
$n>k$. By the fact that
$\frac{d}{dz}Y(e^{\al},z)=\al(z)Y(e^{\al},z)$, we have
$$
e^{\al}(k)\al(0)w=(-k-1)e^{\al}(k)w.$$ Since
$$
[\al(n),e^{\be}(r)]=(\al,\be)e^{\be}(n+r),$$ we have
$$
\al(0)e^{\al}(k)w=(-k-1+(\al,\al))e^{\al}(k)w.$$ Note that
$(\al,\al)\in\Z$, so $-k-1+(\al,\al)\in\Z$. Then (2) follows from
Lemma \ref{l3.2} and the fact that  $e^{\al}(k)w$ generates $W$ by
operators $x(n),e^{\be}(n)$, for $x\in H_{4}, $ $\be\in L$, and
$[d,a]=a$, $[d,b]=-b$. \qed

\begin{theorem}\label{t4.2}
Let $W$ be an irreducible $V_{\widehat{H}_{4}}(1,0)[L]$-module, then
$W\cong \C [L]\otimes U$, where $U$ is an irreducible module of
$V_{\widehat{H}_{4}}(1,0)$ such that for $\al\in L$, $\al(0)$ acts
on $U$ semisimply with only integral eigenvalues and all $\al(n)$
with $n\geq 1$ are locally nilpotent.
\end{theorem}
\pf By Theorem 3.16 of \cite {DLM2} , $W$  viewed as a
$V_{L}$-module is completely reducible and each simple weak module
of $V_{L}$ is $V_{L+\la}$, for some $\la\in L^{\circ}$, where
$L^{\circ}$ is the dual of $L$ and $V_{L}$ is the vertex algebra
associated with the even lattice $L$. Therefore there exist an
element $u\in W$ and $\la\in L^{\circ}$ such that for $\al\in L$,
$$
\al(n)u=0,  \ \al(0)u=(\al,\la)u, $$ and
$$
e^{\al}(-(\al,\la)-1)u\neq 0, \ e^{\al}(-(\al,\la)-1+n)u=0,
$$for $n\in\Z_{+}.$
Since $$ [\al(n),e^{\be}(r)]=(\al,\be)e^{\be}(n+r),$$ it follows
that
$$\al(n)e^{\al}(-(\al,\la)-1)u=0$$
for $n\in\Z_{+}$ and
$$\al(0)e^{\al}(-(\al,\la)-1)u=(\al,\al+\la)e^{\al}(-(\al,\la)-1)u.$$
Let $U$ be the $V_{\H}(1,0)$-module generated by $u$. Then
$\C[L]\otimes U$ is a module of $V_{\H}(1,0)[L]$. For $\al\in L$,
let
$$
E_{\al}=E^{-}(\al,z)Y_{W}(e^{\al},z)E^{+}(\al,z)(-z)^{-\al(0)}. $$
Then as the proof of Theorem 3.10 in \cite{LW}, we have
$$
\frac{d}{dz}E_{\al}=0,$$ and
$$
E_{\al}E_{\be}=\epsilon(\al,\be)E_{\al+\be}.$$ Then it is easy to
see that $W$ is linearly spanned by $\{E_{\al}w|\ \al\in L, w\in
U\}.$

 Define $\psi$:
$\C[L]\otimes U\rightarrow W$ by
$$
\psi(e^{\al}\otimes w)=E_{\al}\otimes w,$$ for $\al\in L,w\in U.$
Quite similar to the proof of Theorem 3.10 in \cite{LW}, we have
$$
\psi(Y_{U[L]}(h(-1){\bf 1},z)(e^{\al}\otimes w))=Y_{W}(h(-1){\bf
1},z)\psi(e^{\al}\otimes w),
$$
$$
\psi(Y_{U[L]}(e^{\al},z)(e^{\be}\otimes
w))=Y_{W}(e^{\al},z)\psi(e^{\al}\otimes w),
$$
for $h\in {\frak h}$, $\al,\be\in L$, $w\in U$. Furthermore, for
$\al=m_{1}c+m_{2}d\in L$, $w\in U$, we have
$$
Y_{U[L]}(a(-1){\bf 1},z)(e^{\al}\otimes w)=z^{m_{2}}e^{\al}\otimes
Y_{U[L]}(a(-1){\bf 1},z)w,$$
$$
Y_{W}(a(-1){\bf 1},z)E_{\al}w=z^{m_{2}}E_{\al}Y_{W}(a(-1){\bf
1},z)w,$$
$$
Y_{U[L]}(b(-1){\bf 1},z)(e^{\al}\otimes w)=z^{-m_{2}}e^{\al}\otimes
Y_{U[L]}(b(-1){\bf 1},z)w,$$
$$
Y_{W}(b(-1){\bf 1},z)E_{\al}w=z^{-m_{2}}E_{\al}Y_{W}(b(-1){\bf
1},z)w.$$ So
$$
\psi(Y_{U[L]}(a(-1){\bf 1},z)(e^{\al}\otimes w))=Y_{W}(a(-1){\bf
1},z)\psi(e^{\al}\otimes w),$$
$$
\psi(Y_{U[L]}(b(-1){\bf 1},z)(e^{\al}\otimes w))=Y_{W}(b(-1){\bf
1},z)\psi(e^{\al}\otimes w).$$ This means that $\psi$ is a module
isomorphism. Hence $W\cong \C[L]\otimes U$.  Since $W$ is
irreducible, it follows that $U$ is irreducible. \qed

\vskip 0.3cm {\bf Acknowledgements} \ We would like to thank the
referee for valuable comments and suggestions.

\end{document}